\documentclass{amsproc}

\usepackage{amsthm}
\usepackage{amssymb}
\usepackage{mathtools}

\usepackage[all]{xy}
\usepackage{enumitem}
\usepackage{color}
\usepackage{url}
\usepackage[normalem]{ulem}

\newcommand{\F}{\mathbb{F}}
\newcommand{\C}{\mathbb{C}}

\newcommand{\G}{\mathbb{G}}
\newcommand{\I}{\mathbb{I}}
\newcommand{\N}{\mathbb{N}}
\newcommand{\Z}{\mathbb{Z}}
\newcommand{\Sym}{\mathbb{S}}

\newcommand{\ad}{\mathrm{ad}}
\newcommand{\id}{\mathrm{id}}

\newcommand{\bm}{\mathbf{m}}
\newcommand{\bp}{\mathfrak{p}}
\newcommand{\bq}{\mathfrak{q}}
\newcommand{\br}{\mathbf{r}}
\newcommand{\bM}{\mathbf{M} }
\newcommand{\nct}{\mathbb{N}_0^\theta}

\newcommand{\pa}{\mathbf p}
\newcommand{\gt}{\widetilde{\mathfrak{g}}}
\newcommand{\g}{\mathfrak{g}}
\newcommand{\h}{\mathfrak{h}}
\newcommand{\n}{\mathfrak{n}}
\newcommand{\ntp}{\widetilde{\mathfrak{n}}_+}
\newcommand{\ntm}{\widetilde{\mathfrak{n}}_-}

\newcommand{\np}{\mathfrak{n}_+}
\newcommand{\nm}{\mathfrak{n}_-}

\newcommand{\rg}{\mathfrak{r}}
\newcommand{\nt}{\widetilde{\mathfrak{n}}}

\newcommand{\cB}{\mathcal{B}}
\newcommand{\cC}{\mathcal{C}}

\newcommand{\cG}{\mathcal{G}}

\newcommand{\cO}{\mathcal{O}}

\newcommand{\cQ}{\mathcal{Q}}
\newcommand{\cR}{\mathcal{R}}
\newcommand{\cU}{\mathcal{U}}

\newcommand{\cW}{\mathcal{W}}
\newcommand{\cX}{\mathcal{X}}

\newcommand{\height}{\operatorname{ht}}
\newcommand{\re}{\operatorname{re}}
\newcommand{\im}{\operatorname{im}}
\newcommand{\GK}{\operatorname{GK-dim}}
\newcommand{\rk}{\operatorname{rank}}

\makeatletter
\numberwithin{equation}{section}
\numberwithin{figure}{section}
\numberwithin{table}{section}

\numberwithin{equation}{section}
\numberwithin{figure}{section}
\theoremstyle{plain}
\newtheorem{thm}{Theorem}[section]
\theoremstyle{plain}

\theoremstyle{plain}

\newtheorem*{conjecture*}{Conjecture}

\newtheorem{example}[thm]{Example}
\theoremstyle{remark}
\newtheorem{rem}[thm]{Remark}

\theoremstyle{plain}
\newtheorem{exa}[thm]{Example}
\newtheorem{algorithm}[thm]{Algorithm}

\makeatother

\begin{document}

\title{Computing finite Weyl groupoids}

\begin{abstract}
In this survey, we present algorithms to compute generalized root systems of Nichols algebras of diagonal type and of contragredient Lie superalgebras. As a consequence, we obtain an algorithm to compute the Lyndon words in the Kharchenko PBW basis associated to each positive root, along with their corresponding hyperwords. These data are essential for obtaining a minimal presentation of Nichols algebras of diagonal type with a finite root system.
\end{abstract}

\author{I. Angiono}
\address[Angiono]{FaMAF-CIEM (CONICET), Universidad Nacional de C\'ordoba, Medina Allende s/n, Ciudad Universitaria, 5000 C\'ordoba, Rep\'ublica Argentina}
\email{ivan.angiono@unc.edu.ar}

\author{L. Vendramin}
\address[Vendramin]{Department of Mathematics and Data Science, Vrije Universiteit Brussel, Pleinlaan 2, 1050 Brussel, Belgium}
\email{Leandro.Vendramin@vub.be}


\maketitle

\section*{Introduction}

Root systems are certain highly symmetric configurations of vectors with rich combinatorics. They play a crucial role in the study of semisimple complex Lie algebras and reductive algebraic groups. The symmetries determine a group acting on the set of vectors, called the Weyl group \cite{Bourbaki}. Through quantizations of these objects, called quantum groups, the root systems also determine important invariants in the theory of this family of Hopf algebras, and the action of the associated Weyl group makes it possible to describe their structure and representations \cite{Lusztig}.

When considering more general types of objects, a more flexible notion of root systems is essential. For example, once we move from semisimple complex Lie algebras to \emph{superalgebras}, the associated Weyl group has many symmetries, but not all: there are so-called \emph{odd} reflections \cite{Serganova-root-system} which also permute the root system, and the notion of a groupoid emerges as a natural way to amalgamate both types of symmetries in a unique structure. 

Several notions of generalized root systems appear in the literature. A good axiomatic framework involving groupoids was introduced in \cite{heckenberger-yamane}, and later refined in \cite{CH-rank3}. There is deep and interesting combinatorics underlying these structures; in rank two (that is, when the vectors belong to $\mathbb{Z}^2$), there exists an infinite family of finite root systems, related to triangulations of polygons \cite{CH-rank2}. When the rank is bigger, there are only finitely many for each rank, and the whole classification was obtained in \cite{CH-classif}. Remarkably, this new notion of generalized root systems has proven to be very useful in understanding Lie superalgebras, first for characteristic zero \cite{heckenberger-yamane} and later on  over arbitrary fields \cite{AA}. 

As there exists a notion of super quantum groups (quantized enveloping superalgebras), it seems natural that a similar behavior appears in the world of Hopf algebras. In this direction, generalized root systems also play a fundamental role in Hopf algebras, more specifically in the study of Nichols algebras \cite{AHS,HS}. Finite root systems lead to the classification of finite-dimensional Nichols algebras over finite groups \cite{Heckenberger-Drinfelddoubles,HV-Classif}, which in turn allowed the classification of finite-dimensional pointed Hopf algebras over complex numbers when the group of group-like elements is abelian \cite{MR2630042,AngGI-Selecta}.

These facts highlight the importance of generalized finite root systems. They are not easy to calculate by hand, and determining the system associated with a given Nichols algebra, or Lie superalgebra, is not easy but is essential on the way to determining a basis and other properties of its structure. In this direction, the present work presents an algorithm that gives a complete description of the root system associated to a Nichols algebra of diagonal type (those appearing for finite abelian groups), or for Lie superalgebras. 

The paper is organized as follows. In Section \ref{sec:Weylgpd} we recall the definition of (generalized) root systems and Weyl groupoids, together with some key results explaining the steps of our algorithm. In 
Section \ref{sec:Nichols+super} we recall the notion of braided vector spaces (with focus on the diagonal type), Nichols algebras, Drinfeld doubles, contragredient Lie superalgebras, and develop the structure of root systems in both contexts. In Section \ref{sec:roots+bases-algorithm} we describe algorithms computing root systems, first for Nichols algebras of diagonal type, and second for contragredient Lie superalgebras; these algorithms are available online. Finally, in Section \ref{sec:examples}, we present two examples of data obtained with the algorithms, one for Nichols algebras and another for Lie superalgebras.

\subsection*{Notation} For each set $\cX$ we denote by $\Sym_{\cX}$ the group of symmetries of $\cX$.
Given a positive integer $N$, $\G_N$ denotes the (multiplicative) subgroup of $N$-th roots of unity of $\C$, and $\G_N'$ is the subset of primitive $N$-th roots of unity.

For each $\theta \in \N$ set $\I_{\theta} \coloneqq \{1,2,\cdots,\theta\}$. If $\theta$ is clear from the context we will simply denote 
$\I = \I_{\theta}$. 
Also, $\{\alpha_i\}_{i\in\I}$ denotes the  standard basis of $\Z^{\I}$. We use multiplicative notation for elements in $\Z_{\geq0}^{\I}$: $\sum_{i\in\I}n_i \alpha_i$ is denoted by $1^{n_1}\cdots \theta^{n_{\theta}}$. For example,
$2\alpha_1+3\alpha_2+\alpha_4 \leftrightsquigarrow 1^22^34$.

Given a super vector space $V=V_0\oplus V_1$ with even component $V_0$ and odd component $V_1$, its superdimension is $\operatorname{sdim} V=(\dim V_0|\dim V_1)$.

\section{Weyl groupoids and root systems}\label{sec:Weylgpd}

We assume that the reader has some familiarity with the theory of groupoids, Coxeter groups and generalized Cartan matrices.

\subsection{Basic data and Coxeter groupoids}\label{subsubsection:basic-data-Coxeter}
Fix $\theta\in\N$ and set $\I = \I_{\theta}$. 

A \emph{basic datum} (of size $\theta$) is a pair $(\cX, \rho)$, where $\cX \neq \emptyset$ is  a  set and $\rho: \I  \to \Sym_{\cX}$ is a function such that $\rho_i^2 = \id$ for all $i\in \I$.

\smallbreak
We associate to a basic datum $(\cX, \rho)$ the quiver
\begin{align*}
\cQ_{\rho} = \{\sigma_i^x := (x, i, \rho_i(x)): i\in \I, x\in \cX\}
\end{align*} 
over $\cX$ (i.e. with set of points $\cX$): each arrow $\sigma_i^x$ has source $s(\sigma_i^x) = \rho_i(x)$ and  target
$t(\sigma_i^x) = x$, $x \in \cX$. 

The diagram of $(\cX, \rho)$ is the graph with points $\cX$ and one edge between $x$ and $y$ decorated with label $i$
for each pair $(x, i, \rho_i(x))$, $(\rho_i(x), i, x)$ such that $x \neq \rho_i(x) = y$. In other words, we replace arrows 
$\sigma_i^x$, $\sigma_i^{\rho_i x}$ for those $x\in\cX$ such that $x\ne\rho_i x$, and omit loops from $\cQ_{\rho}$, that can be deduced from the diagram and $\theta$.
For example, if $\theta = 3$, $\cX=\{x,y\}$, $\rho_1=\rho_2=\id_{\cX}$, $\rho_3(x)=y$, then $\cQ_{\rho}$ and its diagram are
\begin{align}\label{eq:basicdatum-diagram-br3}
&\xymatrix{ \underset{x}{\bullet} \ar@(ur,ul)[]_{1} \ar@(lu,ld)[]_{2} \ar@/^0,5pc/^{3}[r] 
& \ar@/^0,5pc/^{3}[l] \underset{y}{\bullet}\ar@(ur,ul)[]_{1}\ar@(rd,ru)[]_{2} } 
& &\leftrightsquigarrow &
&\xymatrix{ \underset{x}{\circ} \ar@{-}^{3}[r]  & \underset{y}{\circ}.}
\end{align}
We say that $(\cX, \rho)$ is connected when $\cQ_{\rho}$ is connected.

In what follows, we will deal with quotients of the free groupoid $F(\cQ_{\rho})$. In any of these quotients we use the convention
\begin{align*}
\sigma_{i_1}^x\sigma_{i_2}\cdots \sigma_{i_t} =
\sigma_{i_1}^x\sigma_{i_2}^{\rho_{i_1}(x)}\cdots \sigma_{i_t}^{\rho_{i_{t-1}} \cdots \rho_{i_1}(x)}
\end{align*}
i. e., the implicit superscripts are the only ones for which the compositions are defined. Given $i, j\in \I$, $x\in \cX$, $m\in\N$, we introduce a second convention to simplify the notation:
\begin{align*}
(\sigma_i^x\sigma_j)^{m}= 
\underbrace{\sigma_{i}^x\sigma_{j}\cdots \sigma_{i}\sigma_j}_{m \text{ times}}.
\end{align*}

\smallbreak

Keeping in mind the definition of Coxeter groups we introduce an analogous concept for groupoids. 
Recall that:
\begin{itemize} [leftmargin=*]
\item A Coxeter matrix of size $\theta$ is a symmetric matrix $\bm = (m_{ij})_{i,j \in \I}$ with entries in
$\N \cup \{+\infty\}$ such that $m_{ii} = 1$ and $m_{ij} \geq 2$, for all $i\neq j\in \I$.

\item The Coxeter group $\cW(\bm)$ with matrix $\bm$ is the group generated by $\sigma_i$, $i\in\I$, with relations
$(\sigma_i\sigma_j)^{m_{ij}}=e$ for all $i,j\in\I$. In particular, $s_i^2=e$ for all $i\in\I$.
\end{itemize}

Now we move to the groupoid context.
\begin{itemize} [leftmargin=*]
\item A \emph{Coxeter datum} $(\cX, \rho,\bM)$ consists of  a basic datum $(\cX, \rho)$ of size $\theta$ together with a bundle of Coxeter matrices $\bM = (\bm^x)_{x\in \cX}$, $\bm^x = (m^x_{ij})_{i,j\in \I}$, such that
\begin{align}\label{eq:coxeter-datum}
s((\sigma^x_i\sigma_j)^{m^x_{ij}}) &= x, &
m^x_{ij}&=m^{\rho_i(x)}_{ij}, & \text{for all } &x \in \cX, \, i,j \in \I.
\end{align}
\item The \emph{Coxeter groupoid} $\cW = \cW(\cX, \rho, \bM)$ is the groupoid  generated by $\cQ_{\rho}$ with relations
\begin{align}\label{eq:def-coxeter-gpd}
(\sigma_i^x\sigma_j)^{m^x_{ij}} &= \id_x, &i, j\in \I,&  &x\in \cX.
\end{align}
\end{itemize}
Notice that the first equality in \eqref{eq:coxeter-datum} just says that \eqref{eq:def-coxeter-gpd} makes sense since the left-hand side of the equality is a loop based on $x$.
In any Coxeter groupoid, we may speak of the length $\ell$ of a product of $\sigma_i^x$'s and a reduced expression of any element, as we do for Coxeter groups.

\begin{example}\label{ex:coxetergroupoid-base1}
A Coxeter groupoid over a basic datum such that $|\cX|=1$ is simply a Coxeter group.

A Coxeter datum over the basic datum $(\cX, \rho)$ in \eqref{eq:basicdatum-diagram-br3} is determined by two Coxeter matrices 
$\bm^x,\bm^y$ of size 3 such that $m^x_{3j}=m^{y}_{3j}\in 2\Z$ for all $j\in\I$.
\end{example}

\subsection{(Generalized) root systems and Weyl groupoids}\label{subsection:weylgroupoid-def}
As for groups, we will construct groupoids generated by \emph{reflections} which admit a kind of root system 
(a combinatorial discrete datum). In general any one of these groupoids is a Coxeter groupoid, but there are Coxeter groupoids not admitting a root system.

\medbreak

Recall that a \emph{Cartan datum} $\cC = (C^x)_{x\in \cX}$ is a family of generalized Cartan matrices $C^x = (c^x_{ij})_{i,j \in \I}$ such that $c^x_{ij}=c^{\rho_i(x)}_{ij}$ for all $x \in \cX$ and $i,j \in \I$. 

Given $x\in\cX$ and $i\in\I$ we set $s_i^x\in GL (\Z^{\I})$ as the reflection such that
\begin{align}\label{eq:reflection-x}
s_i^x(\alpha_j)&=\alpha_j-c_{ij}^x\alpha_i, & j&\in \I,&  i &\in \I, x \in \cX.
\end{align}
Then the definition of a Cartan datum implies that $s_i^x = s_i^{\rho_i(x)}$.

Let $(\cX, \rho)$ be a connected basic datum of size $\theta$. Recall from  \cite{MR2390080} that a \emph{generalized root system} (GRS) for $(\cX, \rho)$ is a pair $(\cC, \varDelta)$, where 
$\cC = (C^x)_{x\in \cX}$ is a Cartan datum, and 
$\varDelta = (\varDelta^x)_{x\in \cX}$ is a family of subsets $\varDelta^x \subset \Z^{\I}$ (a \emph{bundle of root sets}) such that
\begin{align}
\label{eq:def root system 1}
\varDelta^x &= \varDelta^x_+ \cup \varDelta^x_-, & \varDelta^x_{\pm} &:= \pm(\varDelta^x \cap \Z_{\geq0}^{\I}) \subset \pm\Z_{\geq0}^{\I};
\\ \label{eq:def root system 2}
\varDelta^x \cap \Z \alpha_i &= \{\pm \alpha_i \};& &
\\ \label{eq:def root system 3}
s_i^x(\varDelta^x)&=\varDelta^{\rho_i(x)},& \text{cf. }&\eqref{eq:reflection-x};
\\ \label{eq:def root system 4}
(\rho_i\rho_j)^{m_{ij}^x}(x)&=(x), & m_{ij}^x &:=|\varDelta^x \cap (\Z_{\geq0}\alpha_i+\Z_{\geq0} \alpha_j)|,
\end{align}
for all $x \in \cX$, $i \neq j \in \I$.

We call $\varDelta^x_+$ (resp. $ \varDelta^x_-$)  
the set of \emph{positive} (resp. \emph{negative}) roots.

\medbreak

Given a set $X$ and a group $G$, the set $X \times G  \times X$ admits a structure of groupoid with multiplication defined by
\begin{align*}
(x, g, y)(y,h,z) &= (x, gh, z),& x,y,z &\in X, \, g,h \in G,
\end{align*}
while $(x, g, y)(y',h,z)$ is not defined if $y \neq y'$. 
Using this definition,  the \emph{Weyl groupoid} $\cW(\cR)$  attached to a generalized root system $\cR =  (\cC, \varDelta)$ is
the subgroupoid of $\cX \times GL(\Z^\theta) \times \cX$ generated by 
$\varsigma_i^x = (x, s_i^x,\rho_i(x))$, $i \in \I$, $x \in \cX$.

\medbreak

Let $x, y \in \cX$. Notice that $w(\varDelta^x)= \varDelta^y$ for all $w \in \cW(x, y)$.
Thus it makes sense to consider the sets of \emph{real} and \emph{imaginary} roots at $x$ in analogy with the case of Weyl groups, namely 
\begin{align*}
(\varDelta^{\re})^x &= \bigcup_{y\in \cX}\{ w(\alpha_i): \ i \in \I, \ w \in \cW(y,x) \}, &
(\varDelta^{\im})^x &= \varDelta^{x} - (\varDelta^{\re})^x.
\end{align*}

Notice that a connected basic datum may support several root systems. For example, if $|\cX|=1$, i.e. when we have a Weyl group, and we have an affine Cartan matrix $C$, then we may take $\varDelta$ either as the whole root system of the affine Lie algebra attached to $C$ or else just the subset of real roots. 

This phenomenon does not occur in the case of finite Weyl groupoids, generalizing what happens for Weyl groups: the only possible roots are the real ones.

\begin{thm}\label{th:finite-grs}
\cite[2.11]{CH-rank3}
Let $\cR =  (\cC, \varDelta)$ be a root system. The following statements
are equivalent: 
\begin{enumerate}[leftmargin=*, label=\rm{(\alph*)}]
\item $\cR$ is \emph{finite} (i.e. $\vert \cW \vert < \infty$).
\item $\vert\varDelta^x\vert < \infty$ for all $x\in \cX$.
\item There exists $x\in \cX$ such that $\vert\varDelta^x\vert < \infty$.
\item $\vert(\varDelta^{\re})^x\vert < \infty$ for all $x\in \cX$.
\end{enumerate}
\end{thm}

\begin{example}
Let $(\cX, \rho)$ be the basic datum in \eqref{eq:basicdatum-diagram-br3}. A Cartan datum and a root system $\cR$ for $(\cX, \rho)$ are given by
\begin{align*}
C^x &=\begin{pmatrix} 2 & -1 & 0 \\ -1 & 2 & -1 \\ 0 & -2 & 2 \end{pmatrix}, \qquad C^y =\begin{pmatrix} 2 & -1 & 0 \\ -2 & 2 & -1 \\ 0 & -2 & 2 \end{pmatrix},
\\
\varDelta_{+}^{x}= & \{ 1, 12, 123,  1^22^33^4, 12^23^2, 12^23^3, 12^23^4, 12^33^4, 123^2,  2, 23^2, 23, 3\}, \\
\varDelta_{+}^{y}= & \{ 1, 12^2, 12, 123^2, 12^33^2, 1^22^33^2, 12^23^2, 123, 12^23, 2, 23^2, 23, 3 \}.
\end{align*}
\end{example}

Again by analogy with the case of Weyl/Coxeter groups, there exists a close relationship between Coxeter and Weyl groupoids, which is useful for the determination of reduced expressions on one side, and the computation of (real) roots on the other. 
Given a root system $\cR$ and $x\in \cX$, set $\bm^x = (m^x_{ij})_{i,j\in \I}$, where $m^x_{ij}$ is defined as in \eqref{eq:def root system 4}.
Then $\bM = (\bm^x)_{x\in \cX}$ is a Coxeter datum for $(\cX, \rho)$.
\begin{enumerate}[leftmargin=*, label=\rm{(\roman*)}]
\item By \cite{MR2390080} there is an epimorphism of groupoids $\cW(\cX, \rho, \bM)\to \cW(\cX, \rho, \cC)$. Moreover, if $\cR$ is finite, then this is an isomorphism.

\item \cite[Corollary 3]{MR2390080}
Assume that $\ell(\sigma_{i_1}^x\sigma_{i_2}\cdots \sigma_{i_t})=t$, that is the expression is reduced. For each $j\in\I$,
\begin{align}\label{eq:reduced-expression-extending}
&\ell(\sigma_{i_1}^x\sigma_{i_2}\cdots \sigma_{i_t}\sigma_j)=t+ 1 &&\iff
&&\sigma_{i_1}^x\sigma_{i_2}\cdots \sigma_{i_t}(\alpha_j)\in \varDelta_{+}^{x}.
\end{align}
\end{enumerate}

\smallbreak

Assume now that $\cR$ is finite, and pick $x\in \cX$.

\smallbreak

\begin{enumerate}[leftmargin=*, label=\rm{(\roman*)}, resume]
\item \cite[Corollary 5]{MR2390080}
There is a unique element $\omega_0^x \in \cW$ ending at $x$ of maximal length $\ell$. Moreover, if $y \in \cX$, then $\omega_0^y$ has the same length $\ell$.

\item \cite[Prop. 2.12]{CH-rank3}
Fix a reduced expression $\omega_0^x =\sigma^x_{i_1} \cdots \sigma_{i_\ell}$. Then 
\begin{align}\label{eq:real-roots}
\varDelta^x_+ &= \{\beta_j: j \in \I_\ell \}, &
\text{where }\beta_j &:= s_{i_1}^x\cdots s_{i_{j-1}}(\alpha_{i_j}) \in \varDelta^x, 
j \in \I_\ell. 
\end{align}
Therefore all roots are real and the root system is univocally determined.
\end{enumerate}

\begin{rem}\label{rem:finite-GRS}
The classification of the finite generalized root systems 
was obtained in \cite{CH-classif}. The list contains examples that do not arise from either Nichols algebras or
contragredient Lie superalgebras. This may be viewed as an analogue of the existence of finite Coxeter groups not coming from Lie algebras or related structures.
\end{rem}

\section{Weyl groupoids for Nichols algebras and Lie superalgebras}\label{sec:Nichols+super}

Next we recall the construction of Weyl groupoids and root systems for Nichols algebras of diagonal type and for contragredient Lie superalgebras. We will give some useful information for the implementation of algorithms computing these root systems in Section \ref{sec:roots+bases-algorithm}.

\subsection{The Weyl groupoid of a Nichols algebra}\label{subsec:Weyl-gpd} 

In this subsection, $\Bbbk$ will denote an algebraically closed field of characteristic 0.

A braided vector space $(V,c)$ is a pair where $V$ is a $\Bbbk$-vector space and $c$ is a $\Bbbk$-linear automorphism of $V\otimes V$ which satisfies the braid equation, namely
\begin{align*}
(c\otimes \id_V)(\id_V\otimes c)(c\otimes \id_V)=(\id_V\otimes c)(c\otimes \id_V)(\id_V\otimes c).
\end{align*}

Let $(V,c)$ be a braided vector space. The Nichols algebra $\cB(V)$ of $(V,c)$ is a quotient of the tensor algebra $T(V)$ by an ideal $\mathcal{J}(V)$ generated by $\Z_{\geq0}$-homogeneous elements of degree $\ge 2$: $\mathcal{J}(V)$ is maximal among the Hopf ideals of the braided Hopf algebra $T(V)$, see e.g. \cite{Andrus} for the precise definition and alternative presentations.

\subsubsection{}
Prominent examples of braided vector spaces are those of diagonal type: there exist a basis $(x_i)_{i\in\I}$ and a matrix $\bq=(q_{ij})_{i,j\in\I}$ (called the braiding matrix) with non-zero entries such that
$c(x_i\otimes x_j)=q_{ij}\, x_j \otimes x_i$ for all $i,j\in\I$.
To each braiding matrix $\bq$ we attach a (generalized) Dynkin diagram whose set of vertices is $\I$, each vertex $i$ labeled with the scalar $q_{ii}$, and there exists an edge between $i\ne j\in\I$ iff $\widetilde{q}_{ij}\coloneqq q_{ij}q_{ji}\ne 1$, labeled with this scalar.

We say that $\bq$ is \emph{connected} if the corresponding Dynkin diagram is so. If $\bp=(p_{ij})_{i,j\in\I}$ is another braiding matrix, we say that $\bq$ and $\bp$ are twist-equivalent if their Dynkin diagrams coincide, i.e. if $q_{ii}=p_{ii}$ and $\widetilde{q}_{ij}=\widetilde{p}_{ij}$ for all $i\ne j\in\I$.

\smallbreak

The Nichols algebra of a braided vector space of diagonal type with braiding matrix $\bq$ is simply denoted by $\cB_{\bq}$, and correspondingly $\mathcal{J}(V)$ is denoted by $\mathcal{J}_{\bq}$. The tensor algebra $T(V)$ admits a canonical $\Z^{\I}$-grading of braided Hopf algebras determined by $\deg x_i=\alpha_i$. The braiding preserves this $\Z^{\I}$-grading, so the ideal $\mathcal{J}_{\bq}$ is $\Z^{\I}$-homogeneous. Thus the Nichols algebra $\cB_{\bq}$ inherits the $\Z^{\I}$-grading.

The matrix $\bq$ gives rise to a $\Z$-bilinear form $\bq:\Z^{\I}\times\Z^{\I}\to\Bbbk^{\times}$ such that $\bq(\alpha_i,\alpha_j)=q_{ij}$.
For each $\beta\in\Z^{\I}$ set $\bq_{\beta}\coloneqq \bq(\beta,\beta)$. Given two homogeneous elements $x,y$ in either $T(V)$ or $\cB_{\bq}$, we define the braided bracket as
\begin{align*}
[x,y]_c &\coloneqq xy-\bq(\alpha,\beta)yx, & \alpha=\deg x, \, & \beta=\deg y \in\Z^{\I}.
\end{align*}

\subsubsection{}

Recall that a word in letters $x_i$, $i\in\I$ is a Lyndon word if it is smaller than any of its proper ends. The combinatorics of Lyndon words is very rich and was explored in connection with the determination of bases in several contexts, see e.g. \cite{Lalonde-Ram} for the case of enveloping algebras of Lie algebras.

In the case of Nichols algebras of diagonal type, there exists a way to attach a kind of PBW basis (i.e. a basis made by ordered products) whose letters are labeled by Lyndon words. We recall the main definitions and results following \cite{Kharchenko}. For each Lyndon word $u$ we pick an element $[u]_c\in T(V)$ recursively on the length: if $u=x_i$ then $[u]_c=x_i$, while for $u$ of length $>1$ we set $[u]_c \coloneqq \left[[v]_c,[w]_c \right]_c$, where $u=vw$ is the Shirshov decomposition of $u$. The set of ordered products
\begin{align*}
& [u_1]_c^{n_1}\cdots [u_k]_c^{n_k}, & &k\in\Z_{\geq0}, u_1<\cdots < u_k \text{ Lyndon words}, n_i\in\N,
\end{align*}
is a basis of $T(V)$, so 
the images of these elements under the quotient 
by $\mathcal{J}_{\bq}$) span $\cB_{\bq}$.
According to \cite{Kharchenko} there exists a subset $\cG$ of Lyndon words determined by $\mathcal{J}_{\bq}$ as well as a function 
$h:\cG\to\N\cup\{\infty\}$ such that the set of products
\begin{align*}
& [u_1]_c^{n_1}\cdots [u_k]_c^{n_k}, & &k\in\Z_{\geq0}, u_1<\cdots <  u_k \in\cG, 0<n_i<h(u_i),
\end{align*}
is a basis of $\cB_{\bq}$. Here, for each $u\in\cG$, $h(u)\in \{\operatorname{ord} \bq_{\deg u}, \infty \}$ depends on $\mathcal{J}_{\bq}$.

A PBW basis as above is not unique: for example, different orders of the letters $x_i$ give rise to different PBW generators. Anyway the set  $\varDelta_{+}^{\bq} \subset \Z_{\geq0}^{\I}$ of degrees of PBW generators of a PBW basis does not depend on the chosen PBW basis \cite{MR2207786}, and is called the set of positive roots of $\bq$.

The determination of \emph{good} Lyndon words $\cG$ is a hard problem. We will see a quick recipe to determine them when the associated root system is finite. In this case the height is always finite.

\subsubsection{}

The set $\varDelta_{+}^{\bq}$ is the starting point towards the definition of root systems for Nichols algebras of diagonal type. 
Set $\varDelta_{-}^{\bq}=-\varDelta_{+}^{\bq}$, the negative roots, and $\varDelta^{\bq}=\varDelta_{+}^{\bq}\cup \varDelta_{-}^{\bq}$ the set of all roots of $\bq$.

\smallbreak

Nichols algebras are graded by $\Z_{\geq0}^{\I}$ while $\varDelta^{\bq}\subseteq \Z_{\geq0}^{\I}\cup (-\Z_{\geq0}^{\I})$. Thus, in order to \emph{realize} root systems at the level of algebras we need an extension of $\cB_{\bq}$ with a \emph{negative part}, analogous to the case of semisimple Lie algebras (and their universal enveloping algebras). Let $\cU_{\bq}$ be the Drinfeld double of the bosonization $\cB_{\bq}\#\Bbbk\Z^{\I}$, as considered in \cite{Heckenberger-Drinfelddoubles}. This is a Hopf algebra generated by group-like elements $K_i^{\pm 1}$, $L_i^{\pm 1}$, $i\in\I$, and skew-primitive elements $E_i$, $F_i$, $i\in\I$, which resembles quantized enveloping algebras. Indeed, $\cU_{\bq}$ has a triangular decomposition: the multiplication gives a linear isomorphism $\cU_{\bq}\simeq \cU_{\bq}^+\otimes\cU_{\bq}^0\otimes \cU_{\bq}^-$, where
\begin{itemize}[leftmargin=*]\renewcommand{\labelitemi}{$\diamond$}
\item the subalgebras $\cU_{\bq}^+$, $\cU_{\bq}^0$, and $\cU_{\bq}^-$ are generated by $\{E_i\}_{i\in\I}$, $\{K_i^{\pm1},L_i^{\pm1}\}_{i\in\I}$, and $\{F_i\}_{i\in\I}$, respectively; 
\item $\cU_{\bq}^+\simeq \cB_{\bq}$ and $\cU_{\bq}^-\simeq \cB_{\bq^t}$ as algebras;
\item $\cU_{\bq}^0\simeq \Bbbk[\Z^{\I}\times\Z^{\I}]\simeq\Bbbk[\Z^{2\theta}]$ as Hopf algebra.
\end{itemize}
In addition, $\cU_{\bq}$ is a $\Z^{\I}$-graded Hopf algebra, with
\begin{align*}
\deg E_i&=\alpha_i, & \deg F_i&=-\alpha_i, & \deg K_i&=\deg L_i=0.
\end{align*}

\smallbreak

We say that $\bq$ (or $V$, or   $\cB_{\bq}$) is \emph{admissible}
if for all $i  \neq j$ in $\I$, the set 
$$ \left\{ n \in \Z_{\geq0}: (n+1)_{q_{ii}}(1-q_{ii}^n q_{ij}q_{ji} )=0 \right\}$$
is non-empty. In this case we set
\begin{align}\label{eq:defcij}
	c_{ii}^{\bq} &= 2, & c_{ij}^{\bq}&:= -\min \left\{ n \in \Z_{\geq0}: (n+1)_{q_{ii}}(1-q_{ii}^n q_{ij}q_{ji} )=0 \right\},& i&\neq j.
\end{align}
Then the matrix $C^{\bq} =(c_{ij}^{\bq}) \in \Z^{\I \times \I}$ is a GCM. The definition of integers $c_{ij}^{\bq}$ is related to the associated Nichols algebra:
\begin{align*}
	-c_{ij}^{\bq}&= \max \left\{ n \in \Z_{\geq0}: (\ad_c x_i)^{n}x_j\ne 0 \right\}=\max \{n\in\Z_{\geq0}: \alpha_j+n\alpha_i\in\varDelta_{+}^{\bq}\}.
\end{align*}

Next we consider the reflections $s_i^{\bq}\in GL(\Z^\theta)$ determined by $C^{\bq}$, namely 
\begin{align*}
	s_i^{\bq}(\alpha_j)&=\alpha_j-c_{ij}^{\bq}\alpha_i,& i,j&\in \I.
\end{align*}
To connect them with our algebra $\cU_{\bq}$, for each $i\in \I$ let $\rho_i(\bq)$ be a new braiding matrix given by
\begin{align}\label{eq:rhoiq}
\rho_i(\bq)_{jk} &= \bq(s_i^{\bq}(\alpha_j),s_i^{\bq}(\alpha_k)),& j, k &\in \I.
\end{align}

\begin{thm}\label{th:heck-iso} \cite{Heckenberger-Drinfelddoubles}
There exists an isomorphism of algebras $T_i:\cU_{\rho_i(\bq)} \to \cU_{\bq}$ such that 
$T_i\left((\cU_{\rho_i (\bq)})_{\alpha}\right) =\left(\cU_{\bq}\right)_{s_i^{\bq}(\alpha)}$ for all $\alpha\in\Z^{\I}$.
\end{thm}

Assume now that $\GK \cB_{\bq} < \infty$, or equivalently, that $|\varDelta^{\bq}|<\infty$, by \cite{Angiono-GarciaIglesias}.
Then $\bq$ is admissible by \cite{MR1632802}. By Theorem \ref{th:heck-iso}, $c^\bq_{ij}=c^{\rho_i(\bq)}_{ij}$ for all $i,j \in \I$, and $\GK \cB_{\rho_i(\bq)}=\GK \cB_{\bq} < \infty$ hence $\rho_i(\bq)$ is also admissible. Set
\begin{align*}
\cX_{\bq} &= \{\rho_{i_k} \dots \rho_{i_2}\rho_{i_1}(\bq),& k&\in \Z_{\geq0}, &i_1, \dots, i_k &\in \I\}.
\end{align*}
For each $\bp \in \cX_{\bq}$, $\GK\cB_{\bp}=\GK\cB_{\bq}<\infty$, so $\bp$ is admissible. Thus, the maps in \eqref{eq:rhoiq} restrict to $\rho_i: \cX_{\bq} \to \cX_{\bq}$, and then they satisfy $\rho_i^2 =\id$.

The \emph{generalized root system} of the Nichols algebra $\cB_{\bq}$ is the collection of sets
$(\varDelta^{\bp})_{\bp \in \cX_{\bq}}$. 
The reflections $s_i^{\bp}$ for $i\in \I$ and $\bp \in \cX_{\bq}$ satisfy
\begin{align*}
s_i^{\bp}(\varDelta^{\bp}) &= \varDelta^{\rho_i \bp},
\end{align*}
and we can take the corresponding  Weyl groupoid $\cW_{\bq}$, as a subgroupoid of $\cX_{\bq} \times GL(\Z^{\theta}) \times \cX_{\bq}$ generated by $\sigma_i^{\bp}=(\bp,s_i^{\bp},\rho_i \bp)$, $i\in \I$ and $\bp \in \cX_{\bq}$.

All in all, Nichols algebras of diagonal type with finite dimension or finite GK-dimension give rise to generalized root systems.

\begin{rem}\label{rem:Nichols-finite-GRS}
Nichols algebras with finite root systems were classified in \cite{MR2462836}. The list includes positive parts of quantized enveloping algebras as well as Frobenius-Lusztig kernels, which are examples of the so-called braidings of Cartan type, see Example \ref{ex:Cartan} below. Other examples are those corresponding to Lie superalgebras, either in characteristic zero or in positive characteristic, as summarized in \cite{AA}. The purpose of our algorithms here is to determine whether a braiding matrix $\bq$ is in the list in \cite{MR2462836} \emph{without looking at the list}, and in case the matrix is there, give more information on the structure of $\cB_{\bq}$ when the root system is finite:
\begin{itemize}[leftmargin=*]\renewcommand{\labelitemi}{$\heartsuit$}
\item The computation of positive roots is related to the determination of reduced expressions of the element $w_0^{\bq}$ of maximal length by means of \eqref{eq:reduced-expression-extending} and \eqref{eq:real-roots}. Our algorithm is based on this fact, so it gives the set of positive roots as well as a reduced expression $w_0^{\bq}=\sigma_{i_1}^{\bq}\cdots\sigma_{i_M}$. By \cite{heckenberger-yamane} a set of PBW generators can be obtained using the isomorphisms $T_i$ from Theorem \ref{th:heck-iso}:
\begin{align*}
E_{\beta_k} &:=T_{i_1}^{\bq}\cdots T_{i_{k-1}}(E_{i_k}), & &k=1,\cdots,M.
\end{align*}

\item For the dimension, we need to compute $N_{\beta}=\operatorname{ord} q_{\beta}$. If $N_{\beta}=\infty$ for some $\beta\in\varDelta_{+}^{\bq}$, then $\dim\cB_{\bq}=\infty$; otherwise, $\dim\cB_{\bq}=\prod_{k=1}^{M} N_{\beta_k}$.

\item Once we know the set $\varDelta_{+}^{\bq}$ we can easily get a set of Lyndon words, recursively on the degree. 
As all the roots are real, there exists a unique good Lyndon word $l_{\beta}$ of degree $\beta\in\varDelta_{+}^{\bq}$. 
For $\beta=\alpha_i$, $l_{\alpha_i}=x_i$ as expected. If $\beta\ne \alpha_i$, then there exist $\beta_1,\beta_2\in\varDelta_{+}^{\bq}$ such that $\beta=\beta_1+\beta_2$; by \cite{MR3420518},
\begin{align*}
l_{\beta}=\max \{l_{\beta_1}l_{\beta_2}: \beta_1,\beta_2\in\varDelta_{+}^{\bq}, \beta_1+\beta_2=\beta, l_{\beta_1}<l_{\beta_2}\},
\end{align*}
where we consider the lexicographical order on words. Then we take a decomposition $l_{\beta}=l_{\beta_1}l_{\beta_2}$ as above, with $l_{\beta_1}$ of minimal length, to get $[l_{\beta}]_c=\left[ [l_{\beta_1}]_c, [l_{\beta_2}]_c\right]$.

\item From the classification of finite root systems in \cite{CH-classif}, a finite root system of rank $\theta \ge 3$ has at most $B_{\theta}$ positive roots, where $B_2=14$, $B_3=37$, $B_4=32$, $B_5=49$, $B_6=68$, $B_7=91$, $B_8=120$, and $B_\theta=\theta^2$ for $\theta\ge 9$.
Combining \eqref{eq:reduced-expression-extending} with the previous bound $B_{\theta}$, we can conclude that a braiding matrix $\bq$ of rank $\theta$ admitting an element of the Weyl groupoid of length $>B_{\theta}$ has infinitely many roots.
\end{itemize}

\end{rem}

\begin{exa}\label{ex:Cartan}
A braiding matrix $\bq$ is \emph{of Cartan type} if for all $i\ne j\in\I$ there exists $a_{ij}\in\Z$ such that $q_{ij}q_{ji}=q_{ii}^{a_{ij}}$; if $q_{ii}$ is a root of unity, then we can assume that $-\operatorname{ord} q_{ii}<a_{ij}\le 0$. Set $a_{ii}=2$.
	
If $\bq$ is of Cartan type, then $\rho_i \bq=\bq$ for all $i$, and the Weyl groupoid becomes a group. Moreover, the root system is finite if and only if $A\coloneqq (a_{ij})$ is a finite Cartan matrix \cite{MR1780094}, in which case $\varDelta^{\bq}$ is the root system of $A$. The reflections $\sigma_i^{\bq}$ generate the Weyl group $W$ and isomorphisms $T_i$ on $\cU_{\bq}(\g)$ are those in \cite{Lusztig}.
\end{exa}

Further examples can be found in Section \ref{sec:examples}.

\subsection{The Weyl groupoid of a (modular) Lie (super)algebra}\label{subsec:Weyl-gpd-super}

Next we will recall the definition and properties of the root system of a contragredient Lie superalgebra. As observed by Kac \cite{Kac-super} even over the complex numbers, and later for Lie superalgebras over a field of positive characteristic \cite{Skryabin,BGL}, different matrices (or really, pairs of matrices and parity vectors) may give isomorphic Lie superalgebras, so the Weyl group of Lie algebras over complex numbers is replaced by a Weyl groupoid. The notation and results in this Subsection are mainly from  \cite{AA-super}. We refer to \cite{Musson} for more details on contragredient Lie superalgebras.

\medbreak
Let $\Bbbk$ be a field of characteristic $\ell\ge 0$. Given $A =(a_{ij})\in\Bbbk^{\I\times\I}$ we take a vector space $\mathfrak h$  of dimension $2\theta-\rk A$, with a fixed basis $(h_i)_{ i\in \I_{2\theta-\rk A}}$ and a linearly independent family $(\xi_i)_{i\in \I}$ in $\mathfrak h^*$ such that $\xi_j(h_i)= a_{ij}$, for all $i,j\in\I$.

For each pair $(A,\pa)$, where $\pa=(p_i)\in \G_2^\I$ is the \emph{parity vector}, we define a Lie superalgebra analogous to the Kac-Moody construction\footnote{Indeed we get contragredient Lie algebras when $\pa = (1, \dots, 1)$.}.
First set $\gt:=\gt(A,\pa)$ as the Lie superalgebra generated by $e_i$, $f_i$, $i\in \I$, and $\mathfrak h$,
with relations:
\begin{align*}
[h,h']&=0, &[h,e_i] &= \xi_i(h)e_i, & [h,f_i] &= -\xi_i(h)f_i, & [e_i,f_j]&=\delta_{ij}h_i,
\end{align*}
for all $i, j \in \I$, $h, h'\in \mathfrak h$, where the parity is given by
\begin{align*}
|e_i|&=|f_i|= p_i,&  i&\in\I, & |h|&=1,&  h&\in\mathfrak h.
\end{align*}
This Lie superalgebra is $\Z$-graded, $\gt = \mathop{\oplus}\limits_{k\in \Z}\gt_k$, with $e_i\in \gt_1$, $f_i \in \gt_{-1}$,
$\h = \gt_0$, and this gives a triangular decomposition $\gt=\ntp\oplus\h\oplus\ntm$. The \emph{contragredient Lie superalgebra} associated to $A$, $\pa$ is  
\begin{align*}
\g(A,\pa) \coloneqq \gt(A,\pa) / \rg
\end{align*}
where $\rg=\rg_+\oplus \rg_-$ is the maximal $\Z$-homogeneous ideal intersecting $\h$ trivially \cite{BGL,Musson}.
We identify $e_i$, $f_i$, $h_i$, $\h$ with their images in $\g \coloneqq \g(A,\pa)$.
Clearly $\g$ inherits the grading of $\gt$ and $\g=\np\oplus\h\oplus\nm$, where $\n_\pm =\nt_\pm / \rg_\pm$.
As in \cite{Kac-book,BGL}, we assume from now on that $A$ satisfies
\begin{align}\label{eq:symmetrizable}
&a_{ij}=0\mbox{ if and only if }a_{ji}=0, & \mbox{for all  $i,j\in\I$, $i\neq j$}.
\end{align}
By \cite[Section 4.3]{BGL}, $\g$ is $\Z^{\I}$-graded, with $\deg e_i=-\deg f_i=\alpha_i$, $\deg h=0$, for all $i\in \I$, $h\in \h$.
The set of roots and the sets of positive and negative roots are
\begin{align*}
\nabla^{(A,\pa)} &=\{\alpha\in\Z^{\I} - 0:\g_\alpha\neq0\}, & \nabla_{\hspace{2pt} \pm}^{(A,\pa)} &=\nabla^{(A,\pa)}\cap(\pm\nct).
\end{align*}
As for Lie algebras, there exists an involution exchanging $e_i$ with $f_i$ (sometimes called Chevalley involution), so $\nabla_{\hspace{2pt} -}^{(A,\pa)}= -\nabla_{\hspace{2pt}+}^{(A,\pa)}$. And the existence of a triangular decomposition says that $\nabla^{(A,\pa)}=\nabla_{\hspace{2pt} +}^{(A,\pa)}\cup\nabla_{\hspace{1pt} -}^{(A,\pa)}$.
Also, if the rows of $B$ are obtained from those of $A$ by multiplying them by non-zero scalars, then $\g(A,\pa)\simeq \g(B,\pa)$, so we consider matrices up to this equivalence relation.

\medbreak

As in \cite{Serganova-root-system} we say that $(A,\pa)$ is \emph{admissible} if 
$\ad f_i$ is locally nilpotent in $\g(A,\pa)$ for all  $i\in\I$. For instance, when either $\g(A, \pa)$ is finite-dimensional or $\ell >0$ we get admissible pairs \cite{AA-super}. In this case we can set $C^{(A,\pa)}= \big(c_{ij}^{(A,\pa)}\big)_{i,j\in\I} \in \Z^{\I \times \I}$
as the generalized Cartan matrix with entries $c_{ii}^{(A,\pa)}=2$, and for $i\ne j$,
\begin{align*}
c_{ij}^{(A,\pa)}&:=-\min\{m\in\Z_{\geq0}:(\ad f_i)^{m+1} f_j= 0 \}.
\end{align*}
We have explicit formulas for $c_{ij}^{(A,\pa)}$ depending on $a_{ii}$, $a_{ij}$ and $p_i$. 
To state these formulas, we use the following notation: given $x\in\mathbb{F}_{\ell}=\mathbb Z/\ell \mathbb{Z}$, we set $\widetilde{x}\in \mathbb{Z}$ as the unique representative of $x$ in $\{1-p,2-p,\cdots,-1,0\}$ modulo $\ell$.
For $a_{ii}\ne0$,
\begin{align}\label{eq:cij-super-aii=2}
c_{ij}^{(A,\pa)}&= \begin{cases} \widetilde{a}_{ij}, & a_{ij}\in\F_{\ell}, \text{ either }p_i=1 \text{ or }p_i=-1, \widetilde{a}_{ij} \text{ even},
\\
\widetilde{a}_{ij}-\ell, & a_{ij}\in\F_{\ell}, p_i=-1, \widetilde{a}_{ij} \text{ odd},
\\
1-\tfrac{3-p_i}{2}\ell, & a_{ij}\notin\F_{\ell}, 
\end{cases}
\end{align}
while for $a_{ii}=0$ we have
\begin{align}\label{eq:cij-super-aii=0}
c_{ij}^{(A,\pa)}&= \begin{cases} 0, & a_{ij}=0,
\\
1-\ell, & a_{ij}\ne 0, p_i=1,
\\
-1, & a_{ij}\ne 0, p_i=-1.
\end{cases}
\end{align}
As before, let $s_i^{(A,\pa)}\in GL(\Z^{\I})$ be the involutions given by $C^{(A,\pa)}$.
Following \cite{AA-super} we define for each $i\in \I$ a new pair $\rho_i(A, \pa) := (\rho_iA, \rho_i\pa)$ as follows:
\begin{itemize}
\item $\rho_i\pa = (\overline{p}_j)_{j\in\I}$, $\overline{p}_j = p_jp_i^{c_{ij}^{(A,\pa)}}$. 

\item $\rho_iA =(\overline{a}_{jk})_{j,k\in\I}$, where
\begin{align}\label{eq:aij-super-reflection}
\overline{a}_{jk} &\coloneqq \begin{cases}
c_{ik}^{(A,\pa)}a_{ii}-a_{ik}, & j=i,
\\
c_{ij}^{(A,\pa)}c_{ik}^{(A,\pa)}a_{ji}a_{ii}-c_{ij}^{(A,\pa)}a_{ji}a_{ik}-c_{ik}^{(A,\pa)}a_{ji}a_{ij}+a_{ij}a_{jk}, & j\ne i.
\end{cases}
\end{align}
\end{itemize}

\begin{thm}\label{thm:isomorfismo Ti} 
Let $A\in \Bbbk^{\I\times \I}$ be a matrix satisfying \eqref{eq:symmetrizable} and
$\pa\in (\G_2)^{\I}$ such that $(A,\pa)$ is admissible.
For each $i\in \I$ there exists a Lie superalgebra isomorphism 
$T_i^{(A,\pa)}:\g(\rho_i A,\rho_i\pa)\to\g(A,\pa)$ such that
\begin{align*}
T_i^{(A,\pa)}\left(\g(\rho_i A,\rho_i\pa)_\beta\right) &= \g(A,\pa)_{s_i^{(A,\pa)}(\beta)} &  \text{for all }\beta &\in\pm\Z_{\geq0}^\theta.
\end{align*}
\end{thm}

Let $\cX^{(A, \pa)}$ be the equivalence class of admissible pairs $(A, \pa)$, where $A$ satisfies \eqref{eq:symmetrizable}, with respect to the following equivalence relation $\approx$ 
in $\Bbbk^{\I\times \I}\times \G_2^{\I}$:
\begin{align*}
(A,\pa)&\approx (B,\mathbf{q}) &&\iff && \begin{cases}
\exists k\geq0,i_1,\dots,i_k\in \I,\text{ and a matrix $B'$}\\
\text{such that }\rho_{i_1}\cdots \rho_{i_k}(A,\pa)=(B',\mathbf{q}), \\
\text{and the rows of }B\text{ are obtained from those of }B'\\
\text{by multiplication by nonzero scalars.}
\end{cases}
\end{align*}

Assume that there are reflections $T_i$ for all $(B,\mathbf{q}) \in \cX^{(A, \pa)}$ and  $i\in \I$.
For example, that $\ell>0$ and $(A, \pa)$ satisfy \eqref{eq:symmetrizable}. We set
\begin{align*}
\varDelta^{(B,\mathbf{q})} &= \nabla^{(B,\mathbf{q})} - \{k\, \alpha: \, \alpha\in\nabla^{(B,\mathbf{q})}, k\in\N, k\geq 2\},
\\
\cC^{(A,\pa)} &=\left(\I,\cX^{(A, \pa)},(\rho_i)_{i\in\I},(C^{(B,\mathbf{q})})_{(B,\mathbf{q})\in\cX^{(A, \pa)}}\right).
\end{align*}


\begin{thm}{\cite{AA-super}}\label{thm:root system}
$(\cC^{(A,\pa)}, (\varDelta^{(B,\mathbf{q})})_{(B,\mathbf{q})\in\cX^{(A, \pa)}})$ is a generalized root system. 
\end{thm}

The classification of the finite-dimensional contragredient Lie superalgebras is known either when $\ell =0$ \cite{Kac-super} or when 
$\ell > 0$ \cite{BGL}. The purpose of our algorithms is, in analogy with the case of Nichols algebras, to determine whether \[
\dim\g(A,\pa)<\infty
\]
without looking at the classification, and at the same time to provide more information on the structure of the Lie superalgebra. Among other things:

\begin{itemize}[leftmargin=*]\renewcommand{\labelitemi}{$\heartsuit$}
\item The computation of positive roots. The determination of reduced expressions of the element $w_0^{\bq}$ of maximal length gives the subset $\varDelta_+^{(A,\pa)}$. To compute the whole set of positive roots when $\dim\g(A,\pa)<\infty$ we recall that by \cite{AA},
\begin{align}\label{eq:nabla-g(A,p)}
\nabla_+^{(A,\pa)} = \varDelta_+^{(A,\pa)}\cup\left\{ 2\beta | \beta \in \varDelta_{+,\operatorname{ond}}^{(A,\pa)}\right\}.
\end{align}
where $\varDelta_{+,\operatorname{ond}}^{(A,\pa)}$ is the set of odd positive non-degenerate roots, i.e., those odd positive roots $\beta$ such that the associated non-degenerate invariant bilinear form on $\g(A,\pa)$ does not annihilate $\beta$.
For each $\alpha\in \varDelta_+^{(A,\pa)}$ we may give an expression of $e_{\alpha}$ in terms of the corresponding composition of $T_i$'s as for Nichols algebras, and for $\alpha=2\beta$, $\beta\in\varDelta_{+,\operatorname{ond}}^{(A,\pa)}$, $e_{\alpha}=[e_{\beta},e_{\beta}]$. In addition, odd non-degenerate roots are the images of simple roots in $\varDelta_{+,\operatorname{ond}}^{(A,\pa)}$ under $T_i$'s, so we keep track of which roots are images of these simple roots.

\item For the super dimension, 
\begin{align*}
\dim\g(A,\pa)_0&=2|\nabla_+^{(A,\pa)}|+\dim \h -2|\Delta_{+,\operatorname{odd}}^{(A,\pa)}|, &
\dim\g(A,\pa)_1&=2|\Delta_{+,\operatorname{odd}}^{(A,\pa)}|.
\end{align*}
\item Using the same bound $B_{\theta}$ as for Nichols algebras, we can conclude that a pair $(A,\pa)$ of rank $\theta$ admitting an element of the Weyl groupoid of length $>B_{\theta}$ has infinitely many roots.
\end{itemize}

\section{Computing roots and bases}\label{sec:roots+bases-algorithm}

Next, we describe several algorithms related to generalized root systems. These have been implemented in GAP \cite{GAP} and are freely available in the repository at \url{https://github.com/vendramin/roots}, with DOI \verb+10.5281/zenodo.15554004+.

\subsection{Nichols algebras of diagonal type}
The following algorithms are based on results from \S \ref{subsec:Weyl-gpd}. We allow braiding matrices either over $\C$ or over $\C(t)$.

\begin{algorithm}
If the root system is finite, the algorithm returns the set 
\[
\Delta_+^\bq=\{\beta_1,\dots,\beta_N\}
\]
of positive roots, the set $\mathcal{O}^\bq$ of Cartan roots, the set of
heights $\mathcal{H}^\bq$, $\dim\cB_{\bq}$ and a reduced expression of the longest word
$\omega^\bq$. The convex order $\beta_1<\beta_2<\cdots<\beta_N$ corresponds to 
the reduced expression of $\omega^\bq$, see~\cite{MR3420518}.

Input: A braiding matrix $\bq=(q_{ij})\in\C(t)^{\theta\times\theta}$ 
of diagonal type.

The steps of the algorithm are the following: 
\begin{enumerate}[leftmargin=*]
\item Set $i\leftarrow 1$.
\item Set $M=(M_1|\cdots|M_\theta)\leftarrow \id_{\theta\times\theta}$.
\item Set $\Delta_+^\bq\leftarrow \{\alpha_i\}$.
\item Set $\mathcal{O}^\bq\leftarrow\emptyset$.
\item Set $\mathcal{H}^\bq\leftarrow\emptyset$.
\item Set $\omega^\bq\leftarrow \sigma_1^\bq$\;
\item Repeat the following while at least one column $M_j$ of $M$ belongs to $\Z_{\geq0}^\theta$:
\begin{enumerate}[leftmargin=*]
\item Compute the Cartan matrix $C^\bq=(c^\bq_{jk})$ of $\bq$, where
$c_{jk}^\bq$ is given by \eqref{eq:defcij}. 
\item If $\alpha_i$ is a Cartan root, then add $\omega^\bq(\alpha_i)$ to $\mathcal{O}^\bq$.
\item Compute the reflection $s_i^\bq$.
\item Add the order of $q_{ii}$ to $\mathcal{H}^\bq$ (it could be infinite).
\item $M\leftarrow (s_{i}^\bq M_1|\cdots|s_i^\bq M_\theta)$.
\item $\bq\leftarrow{\br}=(r_{jk})$, where $\br$ is given by \eqref{eq:rhoiq}; namely,
\[
r_{jk}=q_{jk}q_{ik}^{-c_{ij}}q_{ji}^{-c_{ik}}q_{ii}^{c_{ij}c_{ik}}.
\]
\item $i\leftarrow \min\{k\in\{1,\dots,\theta\}\setminus\{i\}\text{: the $k$-th column of $M$ belongs to $\Z_{\geq0}^\theta$}\}$.
\item Add the $i$-th column $M_i$ of $M$ to $\Delta_+^\bq$.
\item Add $\sigma_i$ at the end of the word $\omega^\bq$.
\end{enumerate}
\end{enumerate}

Based on the classification in \cite{CH-classif}, the algorithm terminates if the number of positive roots already computed exceeds the maximum of 250 (the largest size among the exceptional root systems) and the square of the rank.
\end{algorithm}

Once the set of positive roots is computed, various data can be immediately obtained, including the dimension (which is the product of the heights) and the defining relations of the corresponding Nichols algebra.    

The following algorithms compute the Lyndon words and the hyperwords associated
with a set of positive roots. The existence of a PBW basis whose generators are
hyperwords is based on \cite{Kharchenko} and the computation of the Lyndon
words attached to each positive root is based on \cite{MR3420518}, as explained in \S \ref{subsec:Weyl-gpd}.

\begin{algorithm}
\label{alg:lyndon}
This algorithm returns the set $\mathcal{L}$ of Lyndon words associated
with a set $\Delta_+$ of positive roots with convex order. For each $\ell\in\mathcal{L}$ it returns the
set of decompositions of $\ell=\ell_\alpha\ell_\beta$, where 
$\ell_\alpha,\ell_\beta\in\mathcal{L}$ and $\alpha<\beta$ are in $\Delta_+$.

Input: A set $\Delta_+=\{\beta_1,\dots,\beta_N\}$ of positive roots with convex order. 

The steps are the following:
\begin{enumerate}[leftmargin=*]
\item Set $\mathcal{L}\leftarrow\emptyset$.
\item Write $\Delta_+=\{\gamma_1,\dots,\gamma_N\}$, where 
$\height(\gamma_i)\leq\height(\gamma_j)$ if $i<j$.
\item For each $j\in\{1,\dots,N\}$ do the following:
\begin{enumerate}
\item Set $\mathcal{T}\leftarrow\emptyset$. 
\item If $\height(\gamma_j)=1$ then add $x_k$ to $\mathcal{L}$
since $\gamma_j=\alpha_k$ for some
$k\in\{1,\dots,\theta\}$. In this case the set of decompositions of $x_k$ is empty.
\item If $\height(\gamma_j)>1$ then:
\begin{enumerate}
\item $\gamma_j=\beta_k$ for some $k\in\{1,\dots,N\}$. 
\item For each $i\in\{1,\dots,k-1\}$ check whether
$\beta_k-\beta_i\in\Delta_+$. If so, add the pair
$(\beta_i,\beta_t)$ to $\mathcal{T}$, where
$\beta_k-\beta_i=\beta_t$. 
\item Set $\ell_{\beta_k}\leftarrow \max\{\ell_\alpha\ell_\gamma:(\alpha,\gamma)\in\mathcal{T}\}$. Note that
$\ell_\alpha$ and $\ell_\gamma$ were already determined since their heights are $<\height(\beta_k)$. 
\item Add $\ell_{\beta_k}$ to $\mathcal{L}$. The set of decompositions of $\ell_{\beta_k}$ is the
subset of $\mathcal{T}$ consisting of pairs $(\alpha,\gamma)$ such that 
$\ell_{\beta_k}=\ell_\alpha\ell_\gamma$. 
\end{enumerate}
\end{enumerate}
\end{enumerate}
\end{algorithm}

\begin{algorithm}
This algorithm computes the set of hyperwords $\mathcal{K}$ associated with
a set $\Delta_+$ of positive roots with convex order. 

Input: A set $\Delta_+=\{\beta_1,\dots,\beta_N\}$ of positive roots with convex order. 

The steps are the following:
\begin{enumerate}[leftmargin=*]
\item Set $\mathcal{K}\leftarrow\emptyset$. 
\item Use Algorithm~\ref{alg:lyndon} to compute the set $\mathcal{L}$
of Lyndon words and the set of decompositions $\mathcal{D}(\ell)$
of each $\ell\in\mathcal{L}$. 
\item For each $\ell\in\mathcal{L}$ do the following:
\begin{enumerate}
\item If the length of $\ell$ is one, then add $[\ell]_c=x_k$
to $\mathcal{K}$ since $\ell=x_k$ for some
$k\in\{1,\dots,\theta\}$.  
\item If the length of $\ell$ is not one, then $\mathcal{D}(\ell)\ne\emptyset$. 
Thus choose $(\ell_1,\ell_2)\in\mathcal{D}(\ell)$ 
such that the length of $\ell_1$ is minimal and add
$[\ell]_c=[ [\ell_1]_c,[\ell_2]_c]_c$ to $\mathcal{K}$. 
\end{enumerate}
\end{enumerate}
\end{algorithm}

\subsection{Lie superalgebras}

The following algorithms are based on \S \ref{subsec:Weyl-gpd-super}. We can work over fields of arbitrary characteristic.

\begin{algorithm}
If the root system is finite, the algorithm returns the set 
\[
\Delta_+^{(A,\pa)}=\{\beta_1,\dots,\beta_N\}
\]
of positive roots, the set $\Delta_{+,\operatorname{ond}}^{(A,\pa)}$ of odd non-degenerate positive roots and a reduced expression of the longest word $\omega^{(A,\pa)}$.

Input: A pair $(A,\pa)\in\Bbbk^{\theta\times\theta}\times \G_2^{\theta}$ consisting of a matrix and a parity vector.

The steps of the algorithm are the following: 
\begin{enumerate}[leftmargin=*]
	\item Set $i\leftarrow 1$.
	\item Set $M=(M_1|\cdots|M_\theta)\leftarrow \id_{\theta\times\theta}$
	\item Set $\Delta_+^{(A,\pa)}\leftarrow \{\alpha_i\}$.
	\item Set $\Delta_{+,\operatorname{ond}}^{(A,\pa)}\leftarrow\emptyset$.
	\item Set $\omega^{(A,\pa)}\leftarrow \sigma_1^{(A,\pa)}$\;
	\item Repeat the following while at least one column $M_j$ of $M$ belongs to $\Z_{\geq0}^\theta$:
	\begin{enumerate}[leftmargin=*]
		\item Compute the Cartan matrix $C^{(A,\pa)}=(c^{(A,\pa)}_{jk})$ according to \eqref{eq:cij-super-aii=2}, \eqref{eq:cij-super-aii=0}.
		\item If $\alpha_i$ is odd (i.e. $p_i=-1$), then add $\omega^{(A,\pa)}(\alpha_i)$ to $\Delta_{+,\operatorname{odd}}^{(A,\pa)}$; if $\alpha_i$ is also non-degenerate (i.e. $a_{ii}\ne0$), then add $\omega^{(A,\pa)}(\alpha_i)$ to $\Delta_{+,\operatorname{ond}}^{(A,\pa)}$ too.
		\item Compute the reflection $s_i^{(A,\pa)}$.
		\item $M\leftarrow (s_{i}^{(A,\pa)} M_1|\cdots|s_i^{(A,\pa)} M_\theta)$.
		\item ${(A,\pa)}\leftarrow({B}=(b_{jk}),\mathbf{q}=(q_j))$, where $q_j = p_jp_i^{c_{ij}^{(A,\pa)}}$ and $b_{jk}$ is given by \eqref{eq:aij-super-reflection}.
		
		\item $i\leftarrow \min\{k\in\{1,\dots,\theta\}\setminus\{i\}\text{: the $k$-th column of $M$ belongs to $\Z_{\geq0}^\theta$}\}$.
		\item Add the $i$-th column $M_i$ of $M$ to $\Delta_+^{(A,\pa)}$.
		\item Add $\sigma_i$ at the end of the word $\omega^{(A,\pa)}$.
	\end{enumerate}
\end{enumerate}
\end{algorithm}

As before, once the set of positive roots is computed, various data can be immediately obtained.

\section{Examples}\label{sec:examples}

Now we show an example of a Nichols algebra of diagonal type and another of a Lie superalgebra, both of which are finite-dimensional and have the same Weyl groupoid, as expressed in \cite[\S 8.3]{AA}.

\begin{exa}
Let $\xi$ be a primitive cubic root of one and 
\[
\bq=\begin{pmatrix}
        -1&   \xi^2&        1\\ 
         1&  -\xi^2&   \xi^2\\
         1&          1&       -1
 \end{pmatrix},
\]
so $\bq$ is of type $\g(2,3)$. The longest word is 
\[
\omega^{\bq}=\sigma_1\sigma_2\sigma_1\sigma_2\sigma_3\sigma_2\sigma_1\sigma_3\sigma_2\sigma_1.
\]
The set $\Delta_+^{\bq}$ of positive roots consists of 
\begin{align}
\label{eq:g(2,3)d-pos-roots}
& 1, &&12, &&12^2, &&2, &&12^33, && 12^23, &&2^23, && 123, && 23, && 3.
\end{align}
The corresponding vector of heights is $(2,3,2,6,2,3,2,2,3,2)$, so 
\begin{align*}
	\dim \cB_{\bq} = 2^63^36=2^73^4= 10368.
\end{align*}
The set $\cO^{\bq}$ of Cartan roots is 
\[
\{12, 2, 12^23, 23\}.
\]

\end{exa}

\begin{exa}
Let $\Bbbk$ be a field of characteristic 3, $\pa=(-1,-1,-1)$ and
\[
A=\begin{pmatrix}
	0 &  1 & 0 \\ -2 &  2 &   -2 \\ 0 & 1 & 0
\end{pmatrix},
\]
so $(A,\pa)$ gives the Lie superalgebra $\g(2,3)$ \cite{BGL}. The longest word is 
\[
\omega^{(A,\pa)}=\sigma_1\sigma_2\sigma_1\sigma_2\sigma_3\sigma_2\sigma_1\sigma_3\sigma_2\sigma_1.
\]
The set $\Delta_+^{(A,\pa)}$ is exactly \eqref{eq:g(2,3)d-pos-roots}, and $\Delta_{+,\operatorname{ond}}^{(A,\pa)}=\{2\}$, so
\begin{align*}
\nabla_+^{(A,\pa)}=& \{1, 12, 12^2, 2, 2^2, 12^33, 12^23, 2^23, 123, 23, 3\}.
\end{align*}
Thus $\operatorname{sdim}\g(2,3)=(12|14)$.
\end{exa}

\subsection*{Acknowledgement}
This work was initiated at the Universidad de los Andes, Colombia. We thank César Galindo for his kind hospitality.

\section*{Statements and Declarations}
The work of I. A. was partially supported by Conicet and SeCyT (UNC).
The work of L. V. was partially supported by
the project OZR3762 of Vrije Universiteit Brussel, and by 
FWO Senior Research Project G043326N.

The authors have no conflicts of interest relevant to the content of this article.

The algorithms in Section~\ref{sec:roots+bases-algorithm}, implemented in GAP, are available at 
\url{https://github.com/vendramin/roots}, with DOI  \verb+10.5281/zenodo.15554004+.

\bibliographystyle{abbrv}
\bibliography{refs}

\end{document}